\documentclass[10pt,leqno]{article}
\AtBeginDvi{} 
\usepackage{amssymb}
\usepackage{amscd}
\usepackage{amsmath}
\usepackage{amsthm}
\usepackage{mathrsfs}
\usepackage{graphics}
\usepackage{graphicx}
\def\overset#1#2{{\mathrel{\mathop {{#2}_{}}\limits^{#1}}}}
\def\underset#1#2{{\mathrel{\mathop {{}_{} {#2}}\limits_{{#1}_{}}}}}
\def\upplim_#1{\underset{#1}{\overline\lim}\;}
\def\lowlim_#1{\underset{#1}{\underline\lim}\;}
\def\ang#1{{\langle}#1{\rangle}}
\newcommand{\C}{{\mathbf{C}}}

\newcommand{\delbar}{\bar{\partial}}

\newcommand{\sI}{\mathscr{I}}

\newcommand{\lto}{\longrightarrow}

\renewcommand{\O}{{\mathcal{O}}}

\newcommand{\pnc}{{\mathbf{P}^n(\mathbf{C})}}

\usepackage{titlesec}
\titleformat{\section}[hang]
\large
{\bf\thesection.}{.5em}{\bf }[]
\titlespacing{\section}
{0pt}{1.5ex plus .1ex minus .2ex}{0pt}
\numberwithin{equation}{section}

\title{
A Brief Chronicle of the Levi (Hartogs' Inverse) Problem, \\
 Coherence and an Open Problem
}
\date{The University of Tokyo}
\date{(July 2018)}
\author{
Junjiro Noguchi\thanks{
Research supported in part by Grant-in-Aid
 for Scientific Research (C) 15K04917. \hfill\break
\hbox{\quad} AMC2010: 32A99; 32E30 \hfill\break
\hbox{\quad} Key words: coherence, Oka, Levi problem, several complex variables
}
}
\begin{document}
\AtBeginDvi{} 
%
\setlength{\baselineskip}{14.5pt}
\parskip+2.5pt
\maketitle
\thispagestyle{empty}

\begin{abstract}
Here we chronologically summarize briefly the developments
 of the Levi (Hartogs' Inverse) Problem together with the notion
of coherence and its solution,
 shedding light on  some records which have not been
discussed in the past references. In particular, we will discuss
K. Oka's unpublished papers {\em 1943} which solved
 the Levi (Hartogs' Inverse) Problem for unramified Riemann domains
of arbitrary dimension $n \geq 2$, usually referred as
it was solved by Oka IX in {\em 1953}, H.J. Bremermann
and F. Norguet in {\em 1954} for univalent domains, independently.

 At the end we emphasize 
an open problem in a ramified case. 
\end{abstract}

\section{Introduction}
There are now a number of interesting and invaluable comments/surveys
on the developments of the titled Problem and Coherence
in complex analysis of several variables such as, e.g.,
H. Cartan's comments in \cite{oka2}, H. Grauert's Commentary of
 \cite{gr94} Part II, I. Lieb \cite{li}.
The purpose of the present article is to recall briefly the
 developments of the problem and the solution
 together with the notion of coherence,
 shedding light on some records and
unpublished manuscripts of K. Oka that have not been discussed
very much in the former references. 
We will see that the original Levi (Hartogs' Inverse)
 Problem itself was historically solved for unramified
Riemann domains over $\C^n$ in Oka's unpublished papers 1943
(cf.\ \cite{oka43} {\bf Theorem I} at p.\ 27), and then
observe how the notion of
 ``{\em Coherence}'' (\hskip-2pt{\em ``Id\'eaux de domaines ind\'etermin\'es}''
in Oka's terms)  evolved from the problem: Here there is a new point, for
those two issues have been discussed independently in the past
references
(cf., e.g., \cite{grcas} Introduction, \cite{li}).

 As we will see in \S\ref{levi},
the turn of years ``{\em 1943/44}\,'' was indeed a watershed in
the study of analytic function theory of several variables. 
In 1943 K. Oka finished the Three Big Problems in the survey
monograph of Behnke--Thullen 1934  (see items~\ref{bt34},
 \ref{ok43} in \S\ref{levi}). In the next year 1944 K. Oka began to study
the arithmetic property of analytic functions of several variables
by investigating Weierstrass' Preparation Theorem (\cite{okp})
which later led to the notion of ``coherence'' in 1948 (\cite{ok7}),
and in the same year 1944 H. Cartan wrote an experimental paper
\cite{ca44} (cf.\ item~\ref{cart44} in \S\ref{levi}).

We will employ commonly used notion and terminologies in analytic
function theory of several variables  without
definitions (cf., e.g., \cite{guro}, \cite{hor2}, \cite{grr},
 \cite{nis}, \cite{nog16})
except for a {\em Riemann domain} $X$ over $\C^n$
 (resp.\ $\pnc$), which in the present note
is a possibly singular reduced complex space $X$ together with
a holomorphic map $\pi: X \to \C^n$ (resp.\ $\pnc$)
such that the fibers $\pi^{-1}\{z\}$ are discrete
for all $z \in \C^n$ (resp.\ $\pnc$): If $\pi$ is locally biholomorphic,
then $X$ is called an {\em unramified} Riemann domain
over $\C^n$ (resp.\ $\pnc$) ($X$ is necessarily non-singular
in this case).

\section{Levi (Hartogs' Inverse) Problem and coherence}\label{levi}
Karl Weierstrass proved his famous Preparation Theorem about {\em 1860}
(cf.\ \cite{grcas} p.~38). According to K. Oka, K. Weierstrass
considered that the theory of analytic functions of two or more
variables would be quite similar to that of one variable,
 and in particular that the shape of singularities of those functions
 should be arbitrary;
this observation had lasted for quite a while.
Then, however, different phenomena had been found, as the
subject had been studied more.

It is noted that the following list is far from being complete:
\begin{enumerate}
\setlength{\itemsep}{-3pt}
\item
Friedrich Hartogs \cite{har}, {\em 1906}\,:
 He found a phenomenon of simultaneous analytic
     continuation of complex analytic functions of two or more
variables (Hartogs' phenomenon).
\item
Eugenio Elia Levi \cite{lev10}/\cite{lev11}, {\em 1910/11}\,: With the boundary regularity he made
clear the pseudoconvexity property of the boundary of a domain of holomorphy.
\item
Henri Cartan--Peter Thullen \cite{cath}, {\em 1932}\,:
They proved the equivalence of domains of holomorphy and
 holomorphically convex ones.
Then, K. Oka systematically used the property of holomorphic convexity.
\item
Wahlter R\"uckert \cite{ru}, {\em 1933}\,: Here,
R\"uckert's Nullstellensatz, which is sometimes called
the Hilbert--R\"uckert Nullstellensatz, was proved. This result played later
a fundamental role in the study of singular complex analytic spaces and
the coherence, but at the beginning the importance was not recognized
very much.
\item\label{bt34}
Heinrich Behnke--P. Thullen \cite{bt}, {\em 1934}\,:
In this monograph they surveyed the research state of the theory
of several complex variables and raised the {\em Three Big Problems} in
several complex variables, on which they put a special importance:
\begin{enumerate}
\setlength{\itemsep}{-3pt}
\item
Levi (Hartogs' Inverse) Problem\footnote{\hskip1pt
 K. Oka termed the problem as Hartogs' Inverse Problem (cf.\
 \cite{okproc}, \cite{ok6}, \cite{ok9}). Hartogs' Inverse Problem is of
a more primitive or more general form than the Levi Problem in the sense that
the latter assumes a $C^2$ boundary regularity of a given domain,
whereas the first does not.}
(\cite{bt} Chap.\ IV).
\item
Cousin I/II Problems (\cite{bt} Chap.\ V).
\item
Problem of developments (Approximation problem of Runge type)
(\cite{bt} Chap.\ VI).
\end{enumerate}

The monograph was of a special importance for K. Oka to
change his research direction to these problems\footnote{\hskip1pt
  At that time he had been writing an unpublished paper
titled ``Fonctions alg\'ebriques permutables avec une fonction rationnelle
non-lin\'eaire'', pp.\ 97, which was typed in French (cf.\ \cite{okp}).\\
It is also interesting to note that K. Oka solved these problems
in the reversed order in time.}.
\item
Kiyoshi Oka I---III \cite{ok1}, \cite{ok2}, \cite{ok3}, {\em
     1936--1939}:
He solved Problems (b) and (c) above, introducing a  principle
     (method) termed ``J\^oku-Ik\^o''\footnote{\hskip1pt This
 is a method or a principle
of K. Oka all through his series of papers \cite{ok1}---\cite{ok9}
such that to solve a problem on a difficult domain one
embeds the domain into a higher dimensional polydisk, 
extends the problem on the polydisk, and then solves it by making use
of the simple shape of the polydisk (cf.\ \cite{nog16}).}.
The well-known ``Oka Principle'' is in Oka III.
\item
Henri Cartan \cite{ca40}, {\em 1940}\,:
H. Cartan introduced the algebraic notion of ideals, congruence, etc.\
into the theory of analytic function theory, and proved his
matrix decomposition theorem.\footnote{\hskip1pt
This paper was the last one for the research communications that
K. Oka had until the end of the war (1940--'45).}
\item
K. Oka \cite{okproc}, {\em 1941}\,:
This is an announcement of the affirmative solution to the Levi (Hartogs'
     Inverse) Problem for univalent domains (subdomains) of $\C^2$.
\item
K. Oka VI \cite{ok6}, {\em 1942}\,:
This is the full paper of the former one with a remark on the
validity of the result in all dimensions $n \geq 2$.
The key of the proof was the so-called Oka's Heftungslemma which was
proved by means of Weil's integral formula in two dimensional case.
Here he dealt with the Hartogs pseudoconvexity, so that
he put no condition on the regularity of the boundary of the domain,
while the notion of Levi pseudoconvexity needs at least
 $C^2$-regularity.

 In the course of the proof he used a modified
Levi pseudoconvexity by introducing a new class of real-valued functions
called ``fonctions pseudoconvexes'' in Oka VI \S11, which were
 also called ``fonctions plurisousharmoniques''
by Pierre Lelong around the same time.
\item\label{ok43}
K. Oka \cite{okp}, \cite{oka43}, {\em 1943}, Research reports to Teiji Takagi
 (in Japanese, unpublished) (cf.\ \S\ref{ok43-rd}):
 In this year just after Oka VI,
{\em  Oka proved the Levi (Hartogs' Inverse) Problem
 for unramified Riemann domains of
general dimension $\geq 2$} in a series of five research reports
of pp.~109 in total, sent to Teiji Takagi (Tokyo, well known as the
founder of class field theory). The reports were
written in Japanese and unpublished;
 they are now available in \cite{okp} (Japanese).
Let us  quote the most main result of the papers
(see \cite{oka43}\footnote{\hskip1pt For convenience,
 the present author translated the most important last
one among the five into English.}
,  p.\ 27):
\begin{quote}
{\bf Theorem I.}  {\em
A pseudoconvex finite domain\footnote{\hskip1pt Here ``finite domain'' means
``domain over $\C^n$''.}
 with no interior ramification point is a domain of holomorphy.}
\end{quote}

 He remarked this fact three times in his published papers, first
in his survey note \cite{oka49} (1949), in
VIII \cite{ok8} (1951), and in IX \cite{ok9} (1953).

{\em Comparison to Oka VI (1942)}:
The method of the proof was very different from the previous one of Oka VI.
In these reports he proved Heftungslemma by the combination of
J\^oku-Ik\^o and Cauchy's integral formula (in fact, it is
a half of Cauchy integral called the {\em Cousin integral})
 in place of Weil's integral formula, which was not obtained on
an unramified Riemann domain over $\C^n$ (cf.\ \cite{oka43}, \cite{ok9} \S24).

{\em Comparison to Oka IX (1953)}: It is the same in both solutions to
     construct a continuous plurisubharmonic exhaustion
 function on a pseudoconvex
unramified Riemann domain, but in 1943 the coherence theorems
of K. Oka VII/VIII (see item~\ref{ok78}) used in Oka IX was not yet invented,
and hence not used.
It is also noted that in the course he proved a {\em sort of
 coherence theorem} in a special case (cf.\ \S\ref{ok43-rd} below).

{\em N.B.} As H. Cartan wrote in Oka \cite{oka2}, p.~XII,
\vspace{-5pt}
\begin{quote}\sffamily
``Mais, il faut avouer que les aspects techniques de ses d\'emonstrations
 et le mode de pr\'esentation de ses r\'esultats rendent difficile la
t\^ache du lecteur, et que ce n'est qu'au prix d'un r\'eel effort
que l'on parvient \`a saisir la port\'ee de ses r\'esultats,
qui est consid\'erable.''
\end{quote}
\vspace{-5pt}
\begin{quote}\sffamily
``But, one must admit that the technical aspects of his demonstrations
and the mode of presentation of his results make it difficult
for the reader, and that it is  only at the price of a real effort
that one can grasp the extent of his results, which is considerable.''
\textrm{(transl. by the author)}
\end{quote}
it is yet not easy to read these unpublished papers. But, in fact,
it is possible to complete the proofs of the Three Big Problems
without Weierstrass' Preparation Theorem (essential in the proof of
coherence), or the theory of sheaf cohomologies of Cartan--Serre,
nor $L^2-\delbar$ method of H\"ormander (if interested, cf.\ \cite{nog18}).
\item\label{cart44}
H. Cartan \cite{ca44}, {\em 1944}\,:
 Let us quote from Grauert--Remmert \cite{grcas} Introduction 2 ---
\vspace{-5pt}
\begin{quote} \sffamily
Of greatest importance in Complex Analysis is the concept of a
coherent analytic sheaf. Already in 1944 C{\sc artan} had experimented with the
notion of a coherent system of punctual modules. He posed the fundamental
problem, whether for any finite system of holomorphic functions the derived
module system of punctual relations is coherent.
This is exactly the problem, whether the sheaf $\O_{\C^n}$ of germs of holomorphic
functions on complex $n$-space is coherent. In 1948 O{\sc ka}
 gave an affirmative answer; in 1950 C{\sc artan} simplified O{\sc ka}'s proof,
 introducing the terminology "faisceau coh\'erent".
This paper was not known in Japan, in particular to K. O{\sc ka} by
the interruption caused by the war.
\end{quote}
\item\label{hit}
Shin Hitotsumatsu \cite{hi}, {\em 1949}\,:
He generalized Oka's Heftungslemma to the $n$-dimensional case
with arbitrary $n \geq 2$, so that {\em he solved the Levi
(Hartogs' Inverse) Problem in the case of univalent domains of $\C^n$
with $n \geq 2$}.
The proof relied on Weil's integral formula in $n$-variables.
This was published in Japanese and has not been referred in the former references.
\item\label{ok78}
K. Oka VII/VIII \cite{ok7}/\cite{ok8}, {\em 1948\,}\footnote{\hskip1pt
This is the year of the received date of Oka VII \cite{ok7} which was
in fact published in 1950; it took rather long time for publication.
 In a number of references
it is referred so that K. Oka proved his First Coherence Theorem
(the coherence of $\O_{\C^n}$) in {\em 1948} just as in
     item~\ref{cart44}. Here we followed it.
}{\em /51}\,:
He proved his three coherence theorems (1'st, $\O_{\C^n}$; 2'nd, $\sI\ang{A}$,
ideal sheaves of analytic subsets $A$, also proved  by H. Cartan \cite{ca50};
 3'rd, Normalization Theorem).
 Here, {\em Oka's aim of ``coherence'' was to
prove the Levi (Hartogs' Inverse) Problem
obtained in 1943 (item \ref{ok43}) for singular ramified Riemann
domains over $\C^n$}.

This intention of K. Oka, which later countered by
J.E. Forn{\ae}ss' example (see item~\ref{forn}), might have two aspects:
One was the ten years delay of the publication of the solution to the Levi (Hartogs'
     Inverse) Problem, and the other was the motive locomotive of the study
that led to the completely new concept of ``Coherence''
(Id\'eaux de domaines ind\'etermin\'es) in 1948/51.

Nowadays we can find {\em two versions of Oka VII}; one is \cite{ok7},
and the other original is in \cite{oka1}.
The English translation of Oka VII in \cite{oka2} is based on
the original in \cite{oka1}.

Probably, the most notable part of the difference in the two versions is
the last part of the introduction, where in the original \cite{oka1} VII
he wrote:
\vspace{-4pt}
\begin{quote} \sffamily
Or, nous, devant le beau syst\`eme de probl\`emes \`a F. Hartogs et aux 
successeurs, voulons l\'eguer des nouveaux probl\`emes \`a ceux qui nous
suivront; or, comme le champ de fonctions analytiques de plusieurs
variables  s'\'etend heureusement aux divers branches de
math\'ematiques, nous serons permis de r\^ever divers types de nouveaux
probl\`emes y pr\'eparant.
\end{quote}
English translation from \cite{oka2} (only added ``Now'' at the beginning):
\vspace{-4pt}
\begin{quote} \sffamily
Now, having found ourselves face to face with the beautiful problems
 introduced
by F.~H{\sc artogs} and his successors, we should like, in turn, to
 bequeath new problems to those who will follow us. The field of
analytic functions of several variables happily extends into diverse
 branches of mathematics, and we might be permitted to dream of the many
 types of new problems in store for us.\footnote{\hskip1pt
One should notice  that when K. Oka wrote these words, he had
already  finished the Three Big Problems of Behnke--Thullen five years
 before (cf.\ item \ref{ok43}).
}
\end{quote}

 The paragraph above was completely deleted from the published paper
 without notification to K. Oka.
As  K. Oka knew the differences between the published \cite{ok7} and the original
     \cite{oka1} VII, he thought it would be necessary to publish
the original VII, once again in a journal, recognizing that it is of
an extremely exceptional case.
 Thus, he wrote an article, ``{\em Propos post\'erieur}'' \cite{okapp},
which he would have wished to put as an ``Appendix'' to his original VII
(however, in \cite{oka1} VII there is no ``Appendice'').

 It is really interesting to learn how deeply he was concerned with
this problem of the motivation and how he developed his innovative study.
(Cf.\ \cite{nog16} Chap.~9 ``On Coherence'' for more comparisons.)
\item\label{oka9}
K. Oka, IX \cite{ok9}, {\em 1953}\,:
He solved affirmatively the Levi (Hartogs' Inverse) Problem for unramified
Riemann domains over $\C^n$. As mentioned in the paper, the proof was
essentially the same as in his 1943 unpublished papers  (item \ref{ok43}),
and so it was very different from that in Oka VI.
The most essential part,  his Heftungslemma, was here proved
 by making use of his First and Second Coherence Theorems,
J\^oku-Ik\^o, and the Cousin integral (cf.\ item \ref{ok43} and \cite{ok9})
 in place of Weil's integral formula
which was used in Oka VI but not available on unramified Riemann domains
over $\C^n$ (see \cite{ok9} \S23).

 Here he called the problem 
 ``{\em Hartogs' Inverse Problem}''; in fact, he put no condition on the
boundary regularity of the domain.
In the course, he proved (b) and (c) of item~\ref{bt34}
on a holomorphically convex unramified Riemann domain over $\C^n$.

He left two open problems (\cite{ok9} \S23):
 1) the Levi (Hartogs' Inverse) Problem for unramified
Riemann domains over $\pnc$; 2) the Levi (Hartogs' Inverse) Problem for
 {\em ramified} Riemann domains over $\C^n$ or over $\pnc$.
\item\label{brn}
Hans J. Bremermann \cite{br}; Fran\c{c}ois Norguet \cite{nor}, {\em
     1954}\,:
They proved independently the Levi (Hartogs' Inverse) Problem
 for univalent domains of $\C^n$
with arbitrary $n \geq 2$ by proving Oka's Heftungslemma
in the $n$-dimensional case with Weil's integral formula,
as in Hitotsumatsu \cite{hi} 1949 (see item \ref{hit}).
\item
Hans Grauert \cite{gr58}, {\em 1958}\,:
He gave a considerably simplified proof of Oka's Theorem (IX)
 (item \ref{oka9}) by his well-known ``Bumping Method'' combined with
L. Schwartz's finite dimensionality theorem\footnote{\hskip1pt
The proof of this theorem has been known to be rather long and involved 
by making use of the dual spaces.
Now a very simple proof of it in a slightly
generalized form  is available (see \cite{nog16} \S7.3.4.).}.
\item
H. Grauert, about {\em 1960}: A counter-example to the Levi (Hartogs'
     Inverse) Problem for a ramified Riemann domain over $\pnc$.
\item
Reiko Fujita \cite{fuji}, {\em 1963}\,; Akira Takeuchi \cite{tak}, {\em 1964}:
Independently, they affirmatively proved  the Levi (Hartogs'
     Inverse) Problem for unramified Riemann domains over $\pnc$
with at least one boundary point.
\item
Lars H\"ormander \cite{hor}, {\em 1965}\,:
He proved the Levi (Hartogs' Inverse) Problem by solving
directly the $\delbar$-equations with an $L^2$-method.
At the beginning of Chap.\ IV of his well-known book \cite{hor2}
L. H\"ormander wrote: ``In this chapter we abandon the classical methods
     .....
Instead, ..... the Cauchy--Riemann equations where
 the main  point is an $L^2$ estimate .....''

$\displaystyle \vdots$
\item\label{forn}
John Eric Forn{\ae}ss \cite{forn}, {\em 1978}\,:
He gave a counter-example to the Levi (Hartogs' Inverse) Problem
by constructing a smooth 2-sheeted ramified Riemann domain
 $X \overset{\pi}{\lto}\C^2$,
which is locally Stein but not globally: Here being locally Stein
is defined as for every $z \in \C^2$ there is a neighborhood $U$
of $z$ such that $\pi^{-1}U$ is Stein.

$\displaystyle \vdots$
\end{enumerate}

\section{Oka's unpublished papers 1943 and ramified Riemann domains}\label{ok43-rd}

K. Oka solved affirmatively the Levi (Hartogs' Inverse)
Problem for univalent domains of $\C^2$ in 1942 (\cite{ok6}) and for
unramified Riemann
domains over $\C^n$ (arbitrary $n \geq 2$) in 1943  by writing five
 research reports to T. Takagi (Tokyo) (cf.\ \S\ref{levi}, item~\ref{ok43}):
\begin{enumerate}\setlength{\itemsep}{-3pt}
\item[1)]
On analytic functions of several variables: VII -- Subproblem on
congruence of holomorphic functions, pp.~28.
\item[2)]
On analytic functions of several variables: VIII -- The first fundamental
lemma on finite domains without ramification points, pp.~11.
\item[3)]
On analytic functions of several variables: IX -- Pseudoconvex functions,
pp.~30.
\item[4)]
On analytic functions of several variables: X -- The second fundamental
     lemma, pp.~11.
\item[5)]
On analytic functions of several variables: XI -- 
Pseudoconvex domains and  finite domains of holomorphy,
Some theorems on finite domains of holomorphy, pp.~29.
\end{enumerate}

It is noteworthy that in the above VII he proved a special case
of {\em coherence property} (for the so-called Oka maps used in
J\^oku-Ik\^o) already in 1943.

He did not translated these handwritten manuscripts into French for publications,
but immediately began to study the Levi (Hartogs' Inverse) Problem
for Riemann domains with ramifications. He subsequently wrote the
following in the same series as above:
\begin{enumerate}\setlength{\itemsep}{-3pt}
\item[6)]
On analytic functions of several variables: XII -- 
Representation of analytic subsets, pp.~24, 1944.
\item[7)]
On analytic functions of several variables: XII -- 
Extension of the Cousin II problem, pp.~16, 1945.
\item[8)]
XIII -- On a condition in Weierstrass' preparation theorem, pp.~67, 1945.
\end{enumerate}

Here in 6) XII above, he first used Weierstrass' Preparation Theorem for
 the study of the congruence problem of holomorphic
functions:
The purpose was to deal with singular Riemann domains with
ramifications, and this study motivated  and led him to invent the
 ``{\em Coherence''} (\hskip-1.5pt{\em``Id\'eaux de domaines ind\'etermin\'es''}
in Oka's terms)
of holomorphic functions (\cite{ok7}, \cite{ok8}).

In a talk titled ``On analytic functions of several variables''
at Yukawa Institute for Theoretical Physics, Kyoto University 1964,
K. Oka put a special emphasis again on the problem of ramified Riemann
domains (cf.\ \cite{okp}, Unpublished manuscripts, No.\,19), telling that:
\vspace{-4pt}
\begin{quote}  \sffamily
As for Hartogs' Inverse Problem .... And the problem to allow
ramification points remains completely unsolved. I have worked on this
for a rather long time, but I am obstinately keeping the position
to prove it unconditionally. For I have been doing so up to
the present, so otherwise, it is a pity ..... H. Grauert wrote a
paper such that there is an algebraic ramified domain which
 is a domain of holomorphy but not pseudoconvex ....
\end{quote}

H. Grauert also emphasized 
the Levi (Hartogs' Inverse) Problem for ramified Riemann domains
in the talk at Memorial Conference of Kiyoshi Oka's Centennial Birthday,
Kyoto/Nara 2001, which we now formulate as follows.

{\it Problem.} Let $\pi: X \to \C^n$ be a ramified Riemann domain
(cf.\ the beginning of \S1).
 Assume that for every point
$a \in \C^n$ there is a neighborhood $U$ of $a$ in $\C^n$ such
that $\pi^{-}U$ is Stein (locally Stein).
Find sufficient or necessary conditions for $X$ to be Stein.

The problem above is open even for non-singular $X$
(cf.\ \cite{nog17} for some affirmative result).

{
\fontsize{10pt}{0cm}\selectfont

}

\bigskip
\setlength{\baselineskip}{12pt}
\rightline{Graduate School of Mathematical Sciences}
\rightline{University of Tokyo (Emeritus)}
\rightline{Komaba, Meguro-ku, Tokyo 153-8914}
\rightline{Japan}
\rightline{e-mail: noguchi@ms.u-tokyo.ac.jp}
\end{document}